\documentclass[12pt]{amsart}
\usepackage{amscd,amsthm,amsfonts,amssymb,amsmath,euscript}
\usepackage[dvips]{graphicx}
\usepackage[dvips]{graphics}
\usepackage[T2A]{fontenc}
\usepackage[cp866]{inputenc}

\newtheorem{theorem}{Theorem}[section]
\newtheorem{lemma}{Lemma}[section]

\newtheorem{example}{Example}[section]

\newcommand{\Aut}{\mathop{\sf Aut}\nolimits}

\begin{document}

\centerline {\bf Disk single Hurwitz numbers.}

{\ }

\centerline {\bf S.M. Natanzon}

{\ }

\rightline {\it It is devoted to memory of I.M. Gelfand}

{\ }

\section {Introduction.}

The classical Hurwitz numbers are weighted numbers of coverings over
a compact surface with prescribed types of critical values \cite{H}.
For present the full definition let us consider coverings with
ramifications $\varphi: \Omega\rightarrow S$ over compact surface
$S$ of genus $g$ by surfaces $\Omega$. We consider that two
coverings $\varphi$ and $\varphi': \Omega'\rightarrow S$ are
equivalents if there exists a homeomorphism $f:\Omega\rightarrow
\Omega'$ such that $\varphi'f=\varphi$. Denote by $\Aut(\varphi)$
the group of autoequivalence of $\varphi$ and by $|\Aut(\varphi)|$
its order.

The covering $\varphi$ is equivalent to a covering $z\mapsto
z^{n(y)}$ on a neighborhood of any point $y\in\Omega$. Consider
formal variables $\{a_i\}$. The monomial $a_{n(y_1)}...a_{n(y_1)}$,
where $\varphi^{-1}(x)=\{y_1,...,y_k\}$, is called \textit{type of
the value} $x\in S$. A value of a type $a_2a_1^k$ is called
\textit{simple}.

Fix now points $x_1,\ldots,x_v\in S$. Fix also monomials
$a^1,\ldots,a^v$. \textit{The (classical) Hurwitz number} is the sum
$<a^1,..,a^v>_g=\sum\frac{1}{|\Aut(\varphi)|}$, by all equivalent
classes of coverings, that have critical values of types
$a^1,\ldots,a^v$ in the points $x_1,\ldots,x_v$ and have not other
critical values. A Hurwitz number
$<a>^m=<a,a_2a_1^k,..,a_2a_1^k>_0$, where $m+1$ is the number of
critical value is called \textit{classical single Hurwitz number.}
The classical single Hurwitz number connect with the intersection
theory on the moduli space of complex algebraic curves \cite{ELSV}.

Correspond variables $p_i$ to variables $a_i$ and correspond
monomials $p_a=p_{i_1}...p_{i_r}$ to monomials
$a=a_{i_1}...a_{i_r}$. Generating function for classical single
Hurwitz number is
$\Phi(\lambda,p_1,p_2,...)=\sum_{m\geq0}\frac{\lambda^m}{m!}\sum_a
<a>^mp_a$, where second sum is the sum by all monomials. According
to \cite{GJV}, it satisfy to the "cut-and-join" differential
equation
$$\frac{\partial\Phi}{\partial\lambda}=L_{\lambda}\Phi,\;\;\;
L_{\lambda}= \frac{1}{2}\sum_{ij}(i+j) p_ip_j\frac{\partial}{\partial
p_{i+j}}+ \sum_{ij}ij p_{i+j}\frac{\partial^2}{\partial p_i\partial
p_j}.$$

This is a corollary from the fact that classical Hurwitz number
generate a closed topological field theory \cite{D}. The
"cut-and-join" equation has some importent properties. It is in
particulary a differential equation for generating function of Hodge
integrals \cite{Z}.

{\ }

A definition of  Hurwitz numbers for surfaces with boundary and
non-orientable surfaces where suggested in \cite{AN}. In preset
paper I define and investigate \textit{disc single Hurwitz numbers}
They correspond to covering of disk with single non-simple boundary
critical value. It is found, in particulary, recursive equations
defining all disk single Hurwitz numbers.

Generating function $H$ for disc single Hurwitz numbers depend from
complex parameters $\alpha,\beta$ (analog of $\lambda$ from
classical situation) and 4 infinite series of variables
$\acute{p}_i,\grave{p}_i,\bar{p}_i,\dot{p}_i$, describing
topological type of non-simple critical value (analog of $p_i$ from
classical situation). An analog of "cut-and-join" differential
equation $\frac{\partial\Phi}{\partial\lambda}=L_{\lambda}\Phi$ is 2
differential equation, corresponding to interior and boundary simple
critical values.

These differential equation follow, from the fact, that  Hurwitz
numbers for surfaces with boundary form a open-closed
(non-commutative) topological field theory \cite{AN}.

\section {Hurwitz numbers.}

Let $D$ be a closed disk $\{z\in\mathbb{C}||z|\leq 1\}$ with
oriented boundary $\partial D$. We will be consider coverings with
ramification $\varphi: \Omega\rightarrow D$ of degree $k$, where
$\Omega$ is a compact surface with a boundary $\partial \Omega$.

Preimage  $\varphi^{-1}(x)$ of an interior point $x\in D^\circ=
D\setminus\partial D$ consists of $n=n(x)\leq k$ points. Consider a
small simple contour $r\in D^\circ$, around $x$. Its preimage
$\varphi^{-1}(r)$ form simple contours $C_1,...,C_n\in\Omega^\circ$.
The set of degrees of restrictions
$(\deg(\varphi|_{C_1}),...,\deg(\varphi|_{C_n}))$ is called
\textit{(topological) type of interior value $x\in D^\circ$}.
Correspond the monomial $a_1^{t_1}\cdots a_k^{t_k}$ to this
topological type. Here  $a_i$ are commutated formal variables and
$t_i$ is the number of indexes $j$ such that
$\deg(\varphi|_{C_j})=i$ and $a_i^0=1$.

Values of the type $a_1^k$ are called non-critical. All other
interior values are called \textit{interior critical values}. Any
covering has only finite number of interior critical values.
Interior critical values of the type $a_1^{k-1}a_2$ are called
\textit{simple}.

{\ }

A preimage $\varphi^{-1}(y)$ of boundary point $y\in\partial D$ also
consists of $n=n(y)\leq k$ points. Consider a simple small interval
$l\subset D$ around  $y$ with ends $\partial D$. The preimage
$\varphi^{-1}(l)\subset\Omega$ form a graph with $k$ edges. The
vertexes of the graph form two groups that correspond to two ends of
$l$. For convenience, we will call one of the ends and the groups
"left" and other "right", considering that the moving on $\partial
D$ from left to right via $y$ correspond to the orientation of
$\partial D$.

The bipartite graph $\varphi^{-1}(l)$ is called
\textit{(topological) type of boundary critical value $y\in\partial
D$} \cite{AN}. The valency of any vertex of it is not more that 2.
Thus connected components of $\varphi^{-1}(l)$ belong to one of the
next type.
\begin{itemize}
\item $\acute{\textsf{b}}_i$ is a graph with $i$ left and $i+1$ right vertexes;
\item $\grave{\textsf{b}}_i$ is a graph with $i+1$ left and $i$ right vertexes;
\item $\bar{\textsf{b}}_i$ is a graph with $i$ left and $i$ right vertexes;
\item $\dot{\textsf{b}}_i$  is a closed graph with $i$ left and $i$ right vertexes.
\end{itemize}
Correspond to the graphes
$\acute{\textsf{b}}_i,\bar{\textsf{b}}_i,\grave{\textsf{b}}_i,\dot{\textsf{b}}_i$
commutative formal variables
$\acute{b}_i,\bar{b}_i,\grave{b}_i,\dot{b}_i$. Correspond the
monomial to join of graphs. Thus a topological of any boundary value
is a monomial $b=\acute{b}_1^{\acute{s}_1}\cdot\cdot\cdot
\acute{b}_{\acute{n}}^{\acute{s}_{\acute{n}}}
\grave{b}_1^{\grave{s}_1}\cdot\cdot\cdot\grave{b}_{\grave{n}}^{\grave{s}_{\grave{n}}}
\bar{b}_1^{\bar{s}_1}\cdot\cdot\cdot\bar{b}_{\bar{n}}^{\bar{s}_{\bar{n}}}
\dot{b}_1^{\dot{s}_1}\cdot\cdot\cdot\dot{b}_{\dot{n}}^{\dot{s}_{\dot{n}}}$,
where $k=\sum_{i=1}^{\acute{n}}2\acute{s}+$
$\sum_{i=1}^{\bar{n}}(2\bar{s}_i-1)+$
$\sum_{i=1}^{\grave{n}}2\grave{s}_i+$
$\sum_{i=1}^{\dot{n}}2\dot{s}_i$. Denote by $\Aut(b)$ the group of
automorphisms of the graph that correspond to the monomial $b$. Let
$|\Aut(b)|$ be the order of the group.

The changing of the order of the vertexes generates the involution
$b\mapsto b^*$ on the set of monomials. In particulary
$\acute{b_i}^*=\grave{b_i}$, $\bar{b_i}^*=\bar{b_i}$,
$\grave{b_i}^*=\acute{b_i}$, $\dot{b_i}^*=\dot{b_i}$

The values of types $\bar{b}^{\bar{s}}\dot{b}^{\dot{s}}$ are called
non-critical. All other boundary values are called \textit{boundary
critical values}. Any covering has only finite number of boundary
critical values. Boundary critical values of the type
$\acute{b}_1\bar{b}_1^{\bar{n}}\dot{b}_1^{\dot{n}}$ and
$\grave{b}_1\bar{b}_1^{\bar{n}}\dot{b}_1^{\dot{n}}$ are called
\textit{simple boundary critical values}. We call critical values of
types $\acute{b}_1\bar{b}_1^{\bar{n}}\dot{b}_1^{\dot{n}}$ and
$\grave{b}_1\bar{b}_1^{\bar{n}}\dot{b}_1^{\dot{n}}$ as \textit{acute
-points} and \textit{grave-points} respectively.

{\ }

Covering $\varphi_1: \Omega_1\rightarrow D$ ans $\varphi_2:
\Omega_2\rightarrow D$ a called equivalent if there exists a
homeomorphism $\phi: \Omega_1\rightarrow\Omega_2$ such that
$\varphi_1=\varphi_2\phi$. The automorphisms groups
$\Aut(\varphi_i)$ of equivalent coverings are isomorphic. Denote by
$|\Aut(\varphi_i)|$ the order of its.

Fix interior points  $x_1,\ldots,x_v\in D\setminus\partial D$ and
boundary points $y_1,\ldots,y_w\in\partial D$ on a disk $D$. We
consider that the numeration of the boundary points is convenient
with the orientation of $\partial D$. Fix monomials $a^1,\ldots,a^v$
from the variable $a_i$ and  monomials $b^1,\ldots,b^w$ from the
variable $\acute{b}_i,\bar{b}_i,\grave{b}_i,\dot{b}_i$. The number
$<a^1,..,a^v,(b^1,..,b^w)>\sum\frac{1}{|\Aut(\varphi)|}$ is called
\textit{Hurwitz number}. Here the sum is taken by all equivalent
classes of coverings, having critical values in points
$x_1,\ldots,x_v$, $y_1,\ldots,y_w$ of types
$a^1,\ldots,a^v$,$b^1,\ldots,b^w$ and have not other critical
values. The Hurwitz number don't depend from position of critical
points. It is kept by any permutation of $a^i$ cyclic permutation of
$b^i$. Denote by $\{F^{pq}\}$ the the inverse matrix to
$<(\beta_p,\beta_q)>$, where   $\{\beta_p\}$ is the set of all
connected graphs.

\begin{lemma} \label{l1} If $<(c,b)>\neq0$, then $c=b^*$ and
$<(c,b)>= \frac{1}{|\Aut(b)|}$.

If $<(c,\acute{b}_1\bar{b}_1^{\bar{m}}\dot{b}_1^{\dot{m}},b)>\neq0$,
then this is one from the next cases
\begin{itemize}
\item $<(\dot{b}_id^*,\acute{b}_1\bar{b}_1^{\bar{m}}\dot{b}_1^{\dot{m}},
\grave{b}_id)>= \frac{1}{2}\frac{1}{|\Aut(d)|}$;
\item $<(\bar{b}_{i+j}d^*,\acute{b}_1\bar{b}_1^{\bar{m}}\dot{b}_1^{\dot{m}},
\bar{b}_i\grave{b}_jd)>= \frac{1}{|\Aut(d)|}$;
\item $<(\acute{b}_{i+j}d^*,\acute{b}_1\bar{b}_1^{\bar{m}}\dot{b}_1^{\dot{m}},
\grave{b}_i\grave{b}_jd)>=(1-\frac{\delta_{ij}}{2})\frac{1}{|\Aut(d)|}$;
\item $<(\grave{b}_{i+j-1}d^*,\acute{b}_1\bar{b}_1^{\bar{m}}\dot{b}_1^{\dot{m}},
\bar{b}_i\bar{b}_jd)>=
(1-\frac{\delta_{ij}}{2})\frac{1}{|\Aut(d)|}$.
\end{itemize}
If $<(c,\grave{b}_1\bar{b}_1^{\bar{m}}\dot{b}_1^{\dot{m}},b)>\neq0$,
then this is one from the next cases
\begin{itemize}
\item $<(\acute{b}_id,\grave{b}_1\bar{b}_1^{\bar{m}}\dot{b}_1^{\dot{m}},\dot{b}_id^*)>=
\frac{1}{2}\frac{1}{|\Aut(d)|}$;
\item $<(\bar{b}_i\acute{b}_jd,\grave{b}_1\bar{b}_1^{\bar{m}}\dot{b}_1^{\dot{m}},
\bar{b}_{i+j}d^*)>= \frac{1}{|\Aut(d)|}$;
\item $<(\acute{b}_i\acute{b}_jd,\grave{b}_1\bar{b}_1^{\bar{m}}\dot{b}_1^{\dot{m}},
\grave{b}_{i+j}d^*)>=(1-\frac{\delta_{ij}}{2})\frac{1}{|\Aut(d)|};$
\item $<(\bar{b}_i\bar{b}_jd,\grave{b}_1\bar{b}_1^{\bar{m}}\dot{b}_1^{\dot{m}},
\acute{b}_{i+j-1}d^*)>=
(1-\frac{\delta_{ij}}{2})\frac{1}{|\Aut(d)|}$ .
\end{itemize}
\end{lemma}

\textsl{Proof}. The first statement is evident. Let
$<(c,\acute{b}_1\bar{b}_1^{\bar{m}}\dot{b}_1^{\dot{m}},b)>\neq0$.
Then there exists a covering $\varphi: \Omega\rightarrow D$ with
boundary critical values $y_1,y_2,y_3\in\partial D$ of types
$c,\acute{b}_1\bar{b}_1^{\bar{m}}\dot{b}_1^{\dot{m}},b$
(respectively) and without other critical values. Let us consider
points $z_1\in(y_1,y_2)$ and $z_2\in(y_2,y_3)$. The preimage of
$z_1$ consists of $\bar{m}$ simple points, where $\varphi$ is local
homeomorphism, and $\dot{m}+1$ double points where $\varphi$ is
local two-sheeted. The preimage of $z_2$  consists of $\bar{m}+2$
simple points and $\dot{m}$ double points.

Let $p_1,p_2\in\varphi^{-1}(z_2)$ be simple points, corresponding to
$\acute{b}_1$. They come to simple points $q_1,q_2$ of $b$. Then we
have one from the next cases
\begin{itemize}
\item $q_1$ and $q_2$ belong to one connected component of type
$\grave{\textsf{b}}_i$;
\item $q_1$ and $q_2$ belong to connected components of types
$\bar{\textsf{b}}_i$ and $\grave{\textsf{b}}_j$;
\item $q_1$ and $q_2$ belong to different connected components of types
$\bar{\textsf{b}}_i$ and $\bar{\textsf{b}}_j$;
\item $q_1$ and $q_2$ belong to different connected components of types
$\grave{\textsf{b}}_i$ and $\grave{\textsf{b}}_j$.
\end{itemize}

In the first case $b=\grave{b}_id$. Let us consider the restriction
$\varphi': \Omega' D$ of $\varphi$ on the connected component,
containing the points $q_1,q_2$. Consider the restriction
$\varphi'': \Omega'' D$ of $\varphi$ on the complement
$\Omega''=\Omega\setminus\Omega'$. The covering $\varphi''$ has the
critical values $y_1,y_3\in\partial D$ of the types $d^*,d$.
Therefore $\varphi'$ has critical values $y_1,y_2,y_3\in\partial D$
of types $\dot{b}_i,\acute{b}_1e, \grave{b}_i$. Thus $\varphi$ has
critical values $y_1,y_2,y_3\in\partial D$ of types
$\dot{b}_id^*,\acute{b}_1\bar{b}_1^{\bar{m}}\dot{b}_1^{\dot{m}},
\grave{b}_id$.

The equivalent class of $\varphi'$ (respectively $\varphi''$)
contain of all coverings that have the same types of critical values
as $\varphi'$ (respectively $\varphi''$). Moreover,
$\Aut(\varphi)=\Aut(\varphi')\times\Aut(\varphi'')$. Thus
$<(\dot{b}_id^*,\acute{b}_1\bar{b}_1^{\bar{m}}\dot{b}_1^{\dot{m}},
\grave{b}_id)>=$
$\frac{1}{|\Aut(\varphi)|}=\frac{1}{|\Aut(\varphi')|}\frac{1}{|\Aut(\varphi'')|}=$
$\frac{1}{|\Aut(\varphi')|}\frac{1}{|\Aut(d)|}$. The group
$\Aut(\varphi')$ is generated by involution that transpose the
points $p_1$ and $p_2$. Thus
$<(\dot{b}_id^*,\acute{b}_1\bar{b}_1^{\bar{m}}\dot{b}_1^{\dot{m}},
\grave{b}_id)>=$ $\frac{1}{2}\frac{1}{|\Aut(d)|}$.

Analogical arguments prove other cases. The main difference consists
of properties of $\Aut(\varphi')$. This group is non trivial only
for $i=j$ in the last two cases.

Changing of orientation of the boundary of the disk we found that
$<(a,b,c)>=<(c^*,b^*,a^*)>$. Thus the second statement of the lemma
follow from the first.

$\Box$

We denote graphs $\acute{b}_i$ and $\grave{b}_i$ by $\hat{b}_i$ in
situations when differences between them are not important. Denote
by $|b|$ and call length the number of edges in connected graph $b$.
In particulary, $|\bar{b}_i|=2i-1$
$|\hat{b}_i|=|\acute{b}_i|=|\grave{b}_i|=2i$. Denote by $|b|$ the
minimum of lengths of connected components of arbitrary graph $b$.

\begin{lemma} \label{l2}
If $<a_1^ma_2,(c,b)>\neq0$, then this is one from the next cases
\begin{itemize}
\item $<a_1^ma_2,(\dot{b}_i\dot{b}_jd^*,\dot{b}_{i+j}d)>=
(1-\frac{\delta_{ij}}{2})\frac{1}{|\Aut(d)|}$;
\item $<a_1^ma_2,(\dot{b}_i\bar{b}_jd^*,\bar{b}_{i+j}d)>=
\frac{|\bar{b}_j|}{|\Aut(d)|}$;
\item $<a_1^ma_2,(\dot{b}_i\grave{b}_jd^*,\acute{b}_{i+j}d)>=
\frac{|\hat{b}_j|}{|\Aut(d)|}$;
\item $<a_1^ma_2,(\dot{b}_i\acute{b}_jd^*,\grave{b}_{i+j}d)>=
\frac{|\hat{b}_j|}{|\Aut(d)|}$;
\item $<a_1^ma_2,(\bar{b}_i\grave{b}_jd^*,\bar{b}_k\acute{b}_ld)>=
\delta_{(i+j)(k+l)}
\frac{|\bar{b}_i\hat{b}_j\bar{b}_k\hat{b}_l|}{|\Aut(d)|}$;
\item $<a_1^ma_2,(\acute{b}_i\grave{b}_jd^*,\bar{b}_k\bar{b}_ld)>=
\delta_{(i+j+1)(k+l)}
\frac{|\hat{b}_i\hat{b}_j\bar{b}_k\bar{b}_l|}{|\Aut(d)|}$;
\item $<a_1^ma_2,(\grave{b}_i\grave{b}_jd^*,\acute{b}_k\acute{b}_ld)>=
(1-\frac{\delta_{ij}\delta_{kl}}{2})\delta_{(i+j)(k+l)}
\frac{|\hat{b}_i\hat{b}_j\hat{b}_k\hat{b}_l|}{|\Aut(d)|}$;
\item $<a_1^ma_2,(\bar{b}_i\bar{b}_jd^*,\bar{b}_k\bar{b}_ld)>=
(1-\frac{\delta_{ij}\delta_{kl}}{2})\delta_{(i+j)(k+l)}
\frac{|\bar{b}_i\bar{b}_j\bar{b}_k\bar{b}_l|}{|\Aut(d)|},$
\end{itemize}
\end{lemma}

\textsl{Proof}.  Let $<a_1^ma_2,(c,b)>\neq0$. Then there exists a
covering $\varphi: \Omega\rightarrow D$ with boundary critical
values $x,y_1,y_2\in\partial D$ of types $a_1^ma_2,c,b$
(respectively) and without other critical values. Let us consider a
segment $l\subset\partial D$, connecting $y_1$ and $y_2$ by $x$.
Consider also a simple small interval $l_i\subset D$ near $y_i$ with
ends on $\partial D\setminus(\partial D\cap U)$. Then the graph
$b_i=\varphi^{-1}(l_i)$ is the topological type of $y_i$.

Consider a homotopy$\varphi: [1,2]\rightarrow D$  between $l_1$ and
$l_2$ such that $\varphi(t')\cap\varphi(t'')=\emptyset$ for $t'\neq
t''$. Preimage $\varphi^{-1}(x_t)$ of $x_t=\varphi(t)\cap l$
consists of $\deg\varphi$ points if $x_t\neq x$ and  consists
of$\deg\varphi-1$ points if $x_t=x$. This describe the
reconstruction of the graph $b_t=\varphi^{-1}(l_t)$ in the moment,
when $x_t\neq x$. The same reconstruction map the graph $b_1$ to the
graph $b_2^*$. This consists of a cutting of 2 edges of $b_1$ and
the gluing of its with the changing. For realization of this
reconstruction by a covering  $\varphi: \Omega\rightarrow D$ it is
necessary that the orientation of the edges, that appear from order
of vertexes, are generated by an orientation of $\partial D$. All
pears of graphs $b_1,b_2$ of such type and its reconstructions are
presented in lemma \ref{l2}. This gives also and corresponding
Hurwitz numbers.

$\Box$

\section {Disk single Hurwitz numbers.}

We will be say that a boundary critical values $p'$ of a covering
$\varphi$ precedes a boundary critical values $p''$ of $\varphi$ if
the orientation of $\partial D$ generate the orientation of
$(p',p'')$ from $p'$ to $p''$ and $(p',p'')$ does not contain any
critical values.

Denote by $\mathcal{H}(m,\hat{m},b)$ the sen of equivalent classe of
covering with $m$ interior simple critical values, with $\hat{m}$
boundary simple critical values and with single non obviously simple
critical value of type $b$ (that is called \textit{special}). The
set $\mathcal{H}(m,\hat{m},b)$ decompose into subsets
$\mathcal{H}(m,\acute{m},\grave{m},b)$, that consists of coverings
with $m$ interior critical values, $\acute{m}$ acute-points,
$\grave{m}$ grave-points. The set
$\mathcal{H}(m,\acute{m},\grave{m},b)$ decompose into 2 subsets
 $\mathcal{\acute{H}}(m,\acute{m},\grave{m},b)$ and
$\mathcal{\grave{H}}(m,\acute{m},\grave{m},b)$. The first
(respectively second) consists of covering where a acute-point
(respectively a grave-point) precedes the special critical value.

Denote by $\textsc{h}(m,\hat{m},b)$ (respectively
$\textsc{h}(m,\acute{m},\grave{m},b)$,
$\acute{\textsc{h}}(m,\acute{m},\grave{m},b)$,
$\grave{\textsc{h}}(m,\acute{m},\grave{m},b)$) the sum of Hurwitz
numbers from $\mathcal{H}(m,\hat{m},b)$ (respectively
$\mathcal{H}(m,\acute{m},\grave{m},b)$,
$\mathcal{\acute{H}}(m,\acute{m},\grave{m},b)$,
$\mathcal{\grave{H}}(m,\acute{m},\grave{m},b)$). Extend the
definitions of the numbers $\textsc{h}(m,\hat{m},b)$,
$\textsc{h}(m,\acute{m},\grave{m},b)$,
$\acute{\textsc{h}}(m,\acute{m},\grave{m},b)$ and
$\grave{\textsc{h}}(m,\acute{m},\grave{m},b)$ on quotients of
monomials $b=\frac{b^1}{b^2}$. Here we consider that these numbers
are 0, if $b$ is not monomial. The number $\textsc{h}(m,\hat{m},b)$
is called a \textit{disk single Hurwitz number}.

\begin{example} If $\textsc{h}(0,0,b)>0$, then $b=
\bar{b}_1^{\bar{s}_1}\cdot\cdot\cdot\bar{b}_{\bar{n}}^{\bar{s}_{\bar{n}}}
\dot{b}_1^{\dot{s}_1}\cdot\cdot\cdot\dot{b}_{\dot{n}}^{\dot{s}_{\dot{n}}}$
and $\textsc{h}(0,0,b)=$ $\prod_{i=1}^{\bar{n}}\frac{1}{\bar{s}_i!}$
$\prod_{i=1}^{\dot{n}}\frac{1}{i\dot{s}_i!}$.
\end{example}

\begin{lemma} \label{l3.1}
Let
$b=\acute{b}_1^{\acute{s}_1}...\acute{b}_{\acute{n}}^{\acute{s}_{\acute{n}}}
\grave{b}_1^{\grave{s}_1}...\grave{b}_{\grave{n}}^{\grave{s}_{\grave{n}}}
\bar{b}_1^{\bar{s}_1}...\bar{b}_{\bar{n}}^{\bar{s}_{\bar{n}}}
\dot{b}_1^{\dot{s}_1}...\dot{b}_{\dot{n}}^{\dot{s}_{\dot{n}}}.$ Then
$\acute{\textsc{h}}(m,\acute{m}+1,\grave{m},b)=$ $\sum\limits_i
\frac{i}{2}(\dot{s}_i+1)\textsc{h} (m,\acute{m},\grave{m},
b\frac{\dot{b}_i}{\grave{b}_i})+$
$\sum\limits_{ij}((\bar{s}_{i+j}+1)
\textsc{h}(m,\acute{m},\grave{m}, b\frac{\bar{b}_{i+j}}{\bar{b}_i
\grave{b}_j})+$ $(\grave{s}_{i+j}+1)
\textsc{h}(m,\acute{m},\grave{m},
b\frac{\grave{b}_{i+j}}{\grave{b}_i\grave{b}_j})+$
$(\acute{s}_{i+j-1}+1) \textsc{h}(m,\acute{m},\grave{m},
b\frac{\acute{b}_{i+j-1}}{\bar{b}_i\bar{b}_j})).$
\end{lemma}

\textsl{Proof}. Consider a covering $\varphi: \Omega\rightarrow D$
from $\mathcal{\acute{H}}(m,\acute{m}+1,\grave{m},b)$. Denote by $y$
and $y'$ the special point and the precedes of it critical boundary
point respectively. Collapse to a point $y_l$ a segment  $l\in D$
with ends on $\partial D$, that separate the points $y$, $y'$ from
other critical values. Then we have two disks $D'$ and $D''$ and two
coverings $\varphi': \Omega'\rightarrow D'$, $\varphi'':
\Omega''\rightarrow D''$. The critical values of $\varphi'$ are
$y_l$, $y'$ and $y$. Comparing the disks $D$ and $D''$ we see, that
$\varphi''\in\mathcal{H}(m,\acute{m},\grave{m},c)$ for some monomial
$c$. Thus, according to \cite{AN},
$\acute{H}(m,\acute{m}+1,\grave{m},b)=\sum\limits_{pq} H
(m,\acute{m},\grave{m},\beta_p)>$ $F^{pq}<(\beta_q,
\acute{b}_1\bar{b}_1^{\bar{n}}\dot{b}_1^{\dot{n}},b)>$, where
$\{\beta_p\}$  is the set of all bipartite graphes, $\{F^{pq}\}$ is
the invert matrix to $<(\beta_p,\beta_q)>$ and the sum is given by
all pairs of bipartite graphes. It is follow from lemma \ref{l1},
that $F^{pq}=\delta_{\beta_{p},\beta_{q}^*}=|\Aut(\beta_{p})|$.
Moreover it is follow from lemma \ref{l1}, that the sum in the right
parte of decompose to 4 subsumes.  These subsumes depend from forms
of $b$ and $\beta_q=\beta_p^*$.

Let $b=\grave{b}_id$ and $\beta_q=\dot{b}_id^*$. Then
$|\Aut(\beta_p)|=|\Aut(\dot{b}_id)|$
$=\frac{|\Aut(d)|}{i^{\dot{s}_i}\dot{s}_i!}i^{\dot{s}_{i+1}}(\dot{s}_i+1)!=
|\Aut(d)|i(\dot{s}_i+1)$. Moreover
$<(\beta_q,\acute{b}_1\bar{b}_1^{\bar{n}}\dot{b}_1^{\dot{n}},b)>=$
$<(\dot{b}_id^*,\acute{b}_1\bar{b}_1^{\bar{n}}\dot{b}_1^{\dot{n}},\grave{b}_id)>=
\frac{1}{2}\frac{1}{|\Aut(d)|}$ according to lemma \ref{l1}. Thus
this subsume is $\sum\limits_i \frac{i}{2}(\dot{s}_i+1)\textsc{h}
(m,\acute{m},\grave{m}, b\frac{\dot{b}_i}{\grave{b}_i})$ Let
$b=\bar{b}_i\grave{b}_jd$ and $\beta_q=\bar{b}_{i+j}d^*$. Then
$|\Aut(\beta_p)|=|\Aut(\bar{b}_{i+j}d)|$
$=\frac{|\Aut(d)|}{\bar{s}_{i+j}!}(\bar{s}_{i+j}+1)!=
|\Aut(d)|(\bar{s}_{i+j}+1)$. Moreover
$<(\beta_q,\acute{b}_1\bar{b}_1^{\bar{n}}\dot{b}_1^{\dot{n}},b)>=$
$<(\bar{b}_{i+j}d^*,\acute{b}_1\bar{b}_1^{\bar{n}}\dot{b}_1^{\dot{n}},
\bar{b}_i\grave{b}_jd)>=\frac{1}{|\Aut(d)|}$. Thus this subsume is
$\sum\limits_{ij}(\bar{s}_{i+j}+1) \textsc{h}(m,\acute{m},\grave{m},
b\frac{\bar{b}_{i+j}}{\bar{b}_i \grave{b}_j}).$ Let
$b=\grave{b}_i\grave{b}_jd$ and $\beta_q=\acute{b}_{i+j}d^*$. Then
$|\Aut(\beta_p)|=|\Aut(\grave{b}_{i+j}d)|$
$=\frac{|\Aut(d)|}{2^{\grave{s}_{i+j}}\grave{s}_{i+j}!}
2^{(\grave{s}_{i+j}+1)}(\grave{s}_{i+j}+1)!=
2|\Aut(d)|(\grave{s}_{i+j}+1)$. Moreover
$<(\beta_q,\acute{b}_1\bar{b}_1^{\bar{n}}\dot{b}_1^{\dot{n}},b)>=$
$<(\acute{b}_{i+j}d^*,\acute{b}_1\bar{b}_1^{\bar{n}}\dot{b}_1^{\dot{n}},
\grave{b}_i\grave{b}_jd)>=(1-\frac{\delta_{ij}}{2})\frac{1}{|\Aut(d)|}$.
Thus this subsume is $\sum\limits_{ij}(\grave{s}_{i+j}+1)
\textsc{h}(m,\acute{m},\grave{m},
b\frac{\grave{b}_{i+j}}{\grave{b}_i\grave{b}_j}).$ Let
$b=\bar{b}_i\bar{b}_jd$ and $\beta_q=\grave{b}_{i+j-1}d^*$. Then
$|\Aut(\beta_p)|=|\Aut(\acute{b}_{i+j-1}d)|$
$=\frac{|\Aut(d)|}{2^{(\acute{s}_{i+j-1})}\acute{s}_{i+j-1}!}
2^{(\acute{s}_{i+j-1}+1)}(\acute{s}_{i+j-1}+1)!=
2|\Aut(d)|(\acute{s}_{i+j-1}+1)$. Moreover
$<(\beta_q,\acute{b}_1\bar{b}_1^{\bar{n}}\dot{b}_1^{\dot{n}},b)>=$
$<(\grave{b}_{i+j-1}d^*,\acute{b}_1\bar{b}_1^{\bar{n}}\dot{b}_1^{\dot{n}},
\bar{b}_i\bar{b}_jd)>=(1-\frac{\delta_{ij}}{2})\frac{1}{|\Aut(d)|}$.
Thus this subsume is $\sum\limits_{ij}(\acute{s}_{i+j-1}+1)
\textsc{h}(m,\acute{m},\grave{m},
b\frac{\acute{b}_{i+j-1}}{\bar{b}_i\bar{b}_j}).$

$\Box$

\begin{lemma} \label{l3.2}
Let
$b=\acute{b}_1^{\acute{s}_1}...\acute{b}_{\acute{n}}^{\acute{s}_{\acute{n}}}
\grave{b}_1^{\grave{s}_1}...\grave{b}_{\grave{n}}^{\grave{s}_{\grave{n}}}
\bar{b}_1^{\bar{s}_1}...\bar{b}_{\bar{n}}^{\bar{s}_{\bar{n}}}
\dot{b}_1^{\dot{s}_1}...\dot{b}_{\dot{n}}^{\dot{s}_{\dot{n}}}.$ Then
$\grave{\textsc{h}}(m,\acute{m},\grave{m}+1,b)=$
$\sum\limits_{i}(\grave{s}_i+1)\textsc{h}(m,\acute{m},\grave{m},
b\frac{\grave{b}_i}{\dot{b}_i})+$
$\sum\limits_{ij}(2(\bar{s}_i+1)(\grave{s}_j+1)
\textsc{h}(m,\acute{m},\grave{m},b\frac{\bar{b}_i\grave{b}_j}{\bar{b}_{i+j}})+$
$2((\grave{s}_i+1)(\grave{s}_j+1)+ \delta_{ij}(\grave{s}_i+1))
\textsc{h}(m,\acute{m},\grave{m},b\frac{\grave{b}_i\grave{b}_j}{\grave{b}_{i+j}})+$
$\frac{1}{2}((\bar{s}_i+1)(\bar{s}_j+1)+ \delta_{ij}(\bar{s}_i+1))
\textsc{h}(m,\acute{m},\grave{m},b\frac{\bar{b}_i\bar{b}_j}{\acute{b}_{i+j-1}})).$
\end{lemma}

\textsl{Proof}. The proof of this lemma is the same that the
previous , using the second part of lemma \ref{l1}. If $(\beta_q,
\grave{b}_1\bar{b}_1^{\bar{n}}\dot{b}_1^{\dot{n}},b)=$
$(\acute{b}_id^*,\grave{b}_1\bar{b}_1^{\bar{m}}\dot{b}_1^{\dot{m}},\dot{b}_id)$,
Then $|\Aut(\beta_p)|=|\Aut(\grave{b}_id)|=$
$2|\Aut(d)|(\grave{s}_i+1)$ and
$<(\beta_q,\acute{b}_1\bar{b}_1^{\bar{n}}\dot{b}_1^{\dot{n}},b)>=$
$<(\acute{b}_id^*,\grave{b}_1\bar{b}_1^{\bar{m}}\dot{b}_1^{\dot{m}},
\dot{b}_id)>=$ $\frac{1}{2}\frac{1}{|\Aut(d)|}.$ Thus this subsume
is $\sum\limits_{i}(\grave{s}_i+1)\textsc{h}(m,\acute{m},\grave{m},
b\frac{\grave{b}_i}{\dot{b}_i}).$ If  $(\beta_q,
\grave{b}_1\bar{b}_1^{\bar{n}}\dot{b}_1^{\dot{n}},b)=$
$(\bar{b}_i\acute{b}_jd^*,\grave{b}_1\bar{b}_1^{\bar{m}}\dot{b}_1^{\dot{m}},
\bar{b}_{i+j}d)$, then
$|\Aut(\beta_p)|=|\Aut(\bar{b}_i\grave{b}_jd)|=$
$2|\Aut(d)|(\bar{s}_i+1)(\grave{s}_j+1)$ and
$<(\beta_q,\acute{b}_1\bar{b}_1^{\bar{n}}\dot{b}_1^{\dot{n}},b)>=$
$<(\bar{b}_i\acute{b}_jd^*,\grave{b}_1\bar{b}_1^{\bar{m}}\dot{b}_1^{\dot{m}},
\bar{b}_{i+j}d)>= \frac{1}{|\Aut(d)|}.$ Thus this subsume is
$2\sum\limits_{ij}(\bar{s}_i+1)(\grave{s}_j+1)
\textsc{h}(m,\acute{m},\grave{m},b\frac{\bar{b}_i\grave{b}_j}{\bar{b}_{i+j}}).$
If  $(\beta_q,
\grave{b}_1\bar{b}_1^{\bar{n}}\dot{b}_1^{\dot{n}},b)=$
$(\acute{b}_i\acute{b}_jd^*,\grave{b}_1\bar{b}_1^{\bar{m}}\dot{b}_1^{\dot{m}},
\grave{b}_{i+j}d)$, then
$|\Aut(\beta_p)|=|\Aut(\grave{b}_i\grave{b}_jd)|=$
$4|\Aut(d)|((\grave{s}_i+1)(\grave{s}_j+1)+
\delta_{ij}(\grave{s}_i+1))$ and
$<(\acute{b}_i\acute{b}_jd^*,\grave{b}_1\bar{b}_1^{\bar{m}}\dot{b}_1^{\dot{m}},
\grave{b}_{i+j}d)>=(1-\frac{\delta_{ij}}{2})\frac{1}{|\Aut(d)|}$
Thus this subsume is
$2\sum\limits_{ij}((\grave{s}_i+1)(\grave{s}_j+1)+
\delta_{ij}(\grave{s}_i+1))
\textsc{h}(m,\acute{m},\grave{m},b\frac{\grave{b}_i\grave{b}_j}{\grave{b}_{i+j}}).$
If  $(\beta_q,
\grave{b}_1\bar{b}_1^{\bar{n}}\dot{b}_1^{\dot{n}},b)=$
$(\bar{b}_i\bar{b}_jd^*,\grave{b}_1\bar{b}_1^{\bar{m}}\dot{b}_1^{\dot{m}},
\acute{b}_{i+j-1}d)$, then
$|\Aut(\beta_p)|=|\Aut(\bar{b}_i\bar{b}_jd)|=$
$|\Aut(d)|((\bar{s}_i+1)(\bar{s}_j+1)+ \delta_{ij}(\bar{s}_i+1))$
and
$<(\beta_q,\acute{b}_1\bar{b}_1^{\bar{n}}\dot{b}_1^{\dot{n}},b)>=$
$<(\bar{b}_i\bar{b}_jd,\grave{b}_1\bar{b}_1^{\bar{m}}\dot{b}_1^{\dot{m}},
\acute{b}_{i+j-1}d^*)>=
(1-\frac{\delta_{ij}}{2})\frac{1}{|\Aut(d)|}.$ Thus this subsume is
$\frac{1}{2}\sum\limits_{ij}((\bar{s}_i+1)(\bar{s}_j+1)+
\delta_{ij}(\bar{s}_i+1))
\textsc{h}(m,\acute{m},\grave{m},b\frac{\bar{b}_i\bar{b}_j}{\acute{b}_{i+j-1}}).$

$\Box$

The previous lemmas give

\begin{theorem} \label{t3.1} Let
$b=\acute{b}_1^{\acute{s}_1}...\acute{b}_{\acute{n}}^{\acute{s}_{\acute{n}}}
\grave{b}_1^{\grave{s}_1}...\grave{b}_{\grave{n}}^{\grave{s}_{\grave{n}}}
\bar{b}_1^{\bar{s}_1}...\bar{b}_{\bar{n}}^{\bar{s}_{\bar{n}}}
\dot{b}_1^{\dot{s}_1}...\dot{b}_{\dot{n}}^{\dot{s}_{\dot{n}}}.$ Then
$\textsc{h}(m,\hat{m}+1,b)=$ $\sum\limits_i
\frac{i}{2}(\dot{s}_i+1)\textsc{h} (m,\hat{m},
b\frac{\dot{b}_i}{\grave{b}_i})+$
$\sum\limits_{ij}((\bar{s}_{i+j}+1) \textsc{h}(m,\hat{m},
b\frac{\bar{b}_{i+j}}{\bar{b}_i \grave{b}_j})+$ $(\grave{s}_{i+j}+1)
\textsc{h}(m,\hat{m},
b\frac{\grave{b}_{i+j}}{\grave{b}_i\grave{b}_j})+$
$(\acute{s}_{i+j-1}+1) \textsc{h}(m,\hat{m},
b\frac{\acute{b}_{i+j-1}}{\bar{b}_i\bar{b}_j}))+$
$\sum\limits_{i}(\grave{s}_i+1)\textsc{h}(m,\hat{m},
b\frac{\grave{b}_i}{\dot{b}_i})+$
$\sum\limits_{ij}(2(\bar{s}_i+1)(\grave{s}_j+1)
\textsc{h}(m,\hat{m},b\frac{\bar{b}_i\grave{b}_j}{\bar{b}_{i+j}})+$
$2((\grave{s}_i+1)(\grave{s}_j+1)+ \delta_{ij}(\grave{s}_i+1))
\textsc{h}(m,\hat{m},b\frac{\grave{b}_i\grave{b}_j}{\grave{b}_{i+j}})+$
$\frac{1}{2}((\bar{s}_i+1)(\bar{s}_j+1)+ \delta_{ij}(\bar{s}_i+1))
\textsc{h}(m,\hat{m},b\frac{\bar{b}_i\bar{b}_j}{\acute{b}_{i+j-1}})).$
\end{theorem}

Analogical theorem for interior points is

\begin{theorem} \label{t3.2}
Let
$b=\acute{b}_1^{\acute{s}_1}...\acute{b}_{\acute{n}}^{\acute{s}_{\acute{n}}}
\grave{b}_1^{\grave{s}_1}...\grave{b}_{\grave{n}}^{\grave{s}_{\grave{n}}}
\bar{b}_1^{\bar{s}_1}...\bar{b}_{\bar{n}}^{\bar{s}_{\bar{n}}}
\dot{b}_1^{\dot{s}_1}...\dot{b}_{\dot{n}}^{\dot{s}_{\dot{n}}}.$

Then $\textsc{h}(m+1,\hat{m},b)=$

$\sum\limits_{ij}(\frac{i+j}{2}(\dot{s}_{i+j}+1)\textsc{h}(m,\hat{m},b
\frac{\dot{b}_{i+j}}{\dot{b}_i\dot{b}_j})+$
$\frac{ij}{2}((\dot{s}_i+1)
(\dot{s}_j+1)+\delta_{ij}(\dot{s}_i+1))\textsc{h}(m,\hat{m},b
\frac{\dot{b}_i\dot{b}_j}{\dot{b}_{i+j}})+$
$|\bar{b}_j|(\bar{s}_{i+j}+1)\textsc{h}(m,\hat{m},b
\frac{\bar{b}_{i+j}}{\dot{b}_i\bar{b}_j})+$
$i|\bar{b}_j|(\dot{s}_i+1)(\bar{s}_j+1)\textsc{h}(m,\hat{m},b
\frac{\dot{b}_i\bar{b}_j}{\bar{b}_{i+j}})+$
$2|\hat{b}_j|(\acute{s}_{i+j}+1)\textsc{h}(m,\acute{m},\acute{m},b
\frac{\acute{b}_{i+j}}{\dot{b}_i\acute{b}_j})+$
$2i|\hat{b}_j|(\dot{s}_i+1)(\acute{s}_j+1)\textsc{h}(m,\acute{m},\acute{m},b
\frac{\dot{b}_i\acute{b}_j}{\acute{b}_{i+j}})+$
$2|\hat{b}_j|(\grave{s}_{i+j}+1)\textsc{h}(m,\hat{m},b
\frac{\grave{b}_{i+j}}{\dot{b}_i\grave{b}_j})+$
$2i|\hat{b}_j|(\dot{s}_i+1)(\grave{s}_j+1)\textsc{h}(m,\hat{m},b
\frac{\dot{b}_i\grave{b}_j}{\grave{b}_{i+j}}))+$

$\sum\limits_{i+j=k+l}(2|\bar{b}_i\hat{b}_j\bar{b}_k\hat{b}_l|
(\bar{s}_i+1)(\acute{s}_j+1)\textsc{h}(m,\hat{m},b
\frac{\bar{b}_i\acute{b}_j}{\bar{b}_k\acute{b}_l})+$
$2|\bar{b}_i\hat{b}_j\bar{b}_k\hat{b}_l|
(\bar{s}_i+1)(\grave{s}_j+1)\textsc{h}(m,\hat{m},b
\frac{\bar{b}_i\grave{b}_j}{\bar{b}_k\grave{b}_l}))+$

$\sum\limits_{i+j+1=k+l}2(1+\delta_{kl})|\hat{b}_i\hat{b}_j\bar{b}_k\bar{b}_l|
(\acute{s}_i+1)(\grave{s}_j+1)\textsc{h}(m,\hat{m},b
\frac{\acute{b}_i\grave{b}_j}{\bar{b}_k\bar{b}_l})+$

$\sum\limits_{i+j=k+l+1}\frac{1}{2}(1+\delta_{ij})|\bar{b}_i\bar{b}_j\hat{b}_k\hat{b}_l|
((\bar{s}_i+1)(\bar{s}_j+1)+\delta_{ij}(\bar{s}_i+1))
\textsc{h}(m,\hat{m},b
\frac{\bar{b}_i\bar{b}_j}{\acute{b}_k\grave{b}_l})+$

$\sum\limits_{i+j=k+l}(1+\delta_{ij})(1+\delta_{kl})$
$(\frac{1}{4}|\bar{b}_i\bar{b}_j\bar{b}_k\bar{b}_l|
((\bar{s}_i+1)(\bar{s}_j+1)+\delta_{ij}(\bar{s}_i+1))
\textsc{h}(m,\hat{m},b\frac{\bar{b}_i\bar{b}_j}{\bar{b}_k\bar{b}_l})+$

$|\hat{b}_i\hat{b}_j\hat{b}_k\hat{b}_l|
((\acute{s}_i+1)(\acute{s}_j+1)+\delta_{ij}(\acute{s}_i+1))
\textsc{h}(m,\hat{m},b\frac{\acute{b}_i\acute{b}_j}
{\acute{b}_k\acute{b}_l})+$

$|\hat{b}_i\hat{b}_j\hat{b}_k\hat{b}_l|
((\grave{s}_i+1)(\grave{s}_j+1)+\delta_{ij}(\grave{s}_i+1))
\textsc{h}(m,\hat{m},b\frac{\grave{b}_i\grave{b}_j}
{\grave{b}_k\grave{b}_l})),$
\end{theorem}

\textsl{Proof}. Consider a covering $\varphi: \Omega\rightarrow D$
from $\mathcal{\acute{H}}(m+1,\hat{m},b)$. Denote by $y$ and $x$ the
special critical point and an interior critical point respectively.
Collapse to a point $y_l$ a segment  $l\in D$ with ends on $\partial
D$, that separate the points $y$, $x$ from other critical values.
Then we have two disks $D'$ and $D''$ and two coverings $\varphi':
\Omega'\rightarrow D'$, $\varphi'': \Omega''\rightarrow D''$. The
critical values of $\varphi'$ are  $y_l$, $y$ and $x$. Comparing the
disks $D$ and $D''$ we see, that
$\varphi''\in\mathcal{H}(m,\hat{m},c)$ for some monomial $c$. Thus,
according to \cite{AN}, $\acute{H}(m+1,\hat{m},b)=\sum\limits_{pq} H
(m,\hat{m},\beta_p)>$ $F^{pq}<(\beta_q,
\acute{b}_1\bar{b}_1^{\bar{n}}\dot{b}_1^{\dot{n}},b>$.

It is follow from lemma \ref{l2} that the sum from the right parte
decompose on 15 subsumes. Any subsume is calculated similar as in
lemmas \ref{l3.1} and \ref{l3.2}.

$\Box$

Theorems \ref{t3.1}, \ref{t3.2} give an algorithm for calculation of
all Hurwitz numbers $\textsc{h}(m,\hat{m},b)$, starting from
$\textsc{h}(0,0,b)$.

\section {Differential equations on generating the function.}

Consider the algebra of formal power series, generated by
commutative variables $\acute{p}_i,\grave{p}_i,\bar{p}_i,\dot{p}_i$.
Correspond the monomial
$p_b=\acute{p}_1^{\acute{s}_1}...\acute{p}_{\acute{n}}^{\acute{s}_{\acute{n}}}
\grave{p}_1^{\grave{s}_1}...\grave{p}_{\grave{n}}^{\grave{s}_{\grave{n}}}
\bar{p}_1^{\bar{s}_1}...\bar{p}_{\bar{n}}^{\bar{s}_{\bar{n}}}
\dot{p}_1^{\dot{s}_1}...\dot{p}_{\dot{n}}^{\dot{s}_{\dot{n}}}$ to
monomial
$b=\acute{b}_1^{\acute{s}_1}...\acute{b}_{\acute{n}}^{\acute{s}_{\acute{n}}}
\grave{b}_1^{\grave{s}_1}...\grave{b}_{\grave{n}}^{\grave{s}_{\grave{n}}}
\bar{b}_1^{\bar{s}_1}...\bar{b}_{\bar{n}}^{\bar{s}_{\bar{n}}}
\dot{b}_1^{\dot{s}_1}...\dot{b}_{\dot{n}}^{\dot{s}_{\dot{n}}}.$

Consider the generating function
$$\grave{H}(\alpha,\beta,\gamma|
\acute{p}_1,\grave{p}_1,\bar{p}_1,\dot{p}_1,\acute{p}_2,...)=
\sum_{{m,\acute{m}},{\grave{m}}\geq0}
\frac{\alpha^m}{m!}\frac{\beta^{\acute{m}}}{\acute{m}!}
\frac{\gamma^{\grave{m}}}{\grave{m}!} \sum_b
\grave{\textsc{h}}(m,\acute{m},\grave{m},b)p_b,$$ where the second sum is the sum
by all monomial $b$.
In particulary,$\acute{H}(0,0|...)=\grave{H}(0,0|...)=\exp(\bar{p}_1+ \frac{\dot{p}_1}{2})$. Put $H=\acute{H}+\grave{H}$

\begin{theorem}\label{t4.1} $$\frac{\partial \acute{H}}
{\partial\beta}= L_{\beta}H,$$где
$$L_{\beta}= \sum_{i}
i\grave{p}_i\frac{\partial }{\partial\dot{p}_i} +
\sum_{ij}
(\bar{p}_i\grave{p}_j\frac{\partial }{\partial\bar{p}_{i+j}} +
\grave{p}_i\grave{p}_j\frac{\partial }{\partial\grave{p}_{i+j}}+
\bar{p}_i\bar{p}_j\frac{\partial }{\partial\acute{p}_{i+j-1}})$$
\end{theorem}

\textsl{Proof}. Let
$b=\acute{b}_1^{\acute{s}_1}...\acute{b}_{\acute{n}}^{\acute{s}_{\acute{n}}}
\grave{b}_1^{\grave{s}_1}...\grave{b}_{\grave{n}}^{\grave{s}_{\grave{n}}}
\bar{b}_1^{\bar{s}_1}...\bar{b}_{\bar{n}}^{\bar{s}_{\bar{n}}}
\dot{b}_1^{\dot{s}_1}...\dot{b}_{\dot{n}}^{\dot{s}_{\dot{n}}}$,
where $\acute{n}=\acute{n}(b),\acute{s}_i=\acute{s}_i(b)$,
$\grave{n}=\grave{n}(b),\grave{s}_i=\grave{s}_i(b)$
$\bar{n}=\bar{n}(b),\bar{s}_i=\bar{s}_i(b)$
$\dot{n}=\dot{n}(b),\dot{s}_i=\dot{s}_i(b)$.

It is follow from theorem \ref{t3.1} that

$\frac{\partial \acute{H}} {\partial\beta}=$
$\sum_{{m,\acute{m}},{\grave{m}}\geq0}\frac{\alpha^m}{m!}\frac{\beta^{\acute{m}}}
{\acute{m}!}
\frac{\gamma^{\grave{m}}}{\grave{m}!} \sum_b
\acute{\textsc{h}}(m,\acute{m}+1,\grave{m},b)p_b=$

$\sum_{{m,\acute{m},\grave{m}}\geq0}\frac{\alpha^m}{m!}\frac{\beta^{\acute{m}}}
{\acute{m}!}\frac{\gamma^{\grave{m}}}{\grave{m}!}
\sum_b\sum_{i=1}^{\infty} i(\dot{s}_i+1)\textsc{h}
(m,\acute{m},\grave{m}, b\frac{\dot{b}_i}{\grave{b}_i})p_b+$

$\sum_{{m,\acute{m},\grave{m}}\geq0}\frac{\alpha^m}{m!}\frac{\beta^{\acute{m}}}
{\acute{m}!}\frac{\gamma^{\grave{m}}}{\grave{m}!}
\sum_b\sum_{i=1}^{\infty}\sum_{j=1}^{\infty}(\bar{s}_{i+j}+1)
\textsc{h}(m,\acute{m},\grave{m}, b\frac{\bar{b}_{i+j}}{\bar{b}_i
\grave{b}_j})p_b+$

$\sum_{{m,\acute{m},\grave{m}}\geq0}\frac{\alpha^m}{m!}\frac{\beta^{\acute{m}}}
{\acute{m}!}\frac{\gamma^{\grave{m}}}{\grave{m}!}
\sum_b\sum_{i=1}^{\infty}\sum_{j=1}^{\infty}(\grave{s}_{i+j}+1)
\textsc{h}(m,\acute{m},\grave{m},
b\frac{\grave{b}_{i+j}}{\grave{b}_i \grave{b}_j})p_b+$

$\sum_{{m,\acute{m},\grave{m}}\geq0}\frac{\alpha^m}{m!}\frac{\beta^{\acute{m}}}
{\acute{m}!}\frac{\gamma^{\grave{m}}}{\grave{m}!}
\sum_b\sum_{i=1}^{\infty}\sum_{j=1}^{\infty}(\acute{s}_{i+j-1}+1)
\textsc{h}(m,\acute{m},\grave{m},
b\frac{\acute{b}_{i+j-1}}{\bar{b}_i \bar{b}_j})p_b$

Make the changing $c=b\frac{\dot{b}_i}{\grave{b}_i}$ in the first
summand. Then $p_c=p_b\frac{\dot{p}_i}{\grave{p}_i}$ and
$p_b=p_c\frac{\grave{p}_i}{\dot{p}_i}=\grave{p}_i\frac{\partial
p_c}{\partial\dot{p}_i}(\dot{s}_i+1)^{-1}$. Thus
$\sum_b\sum_{i=1}^{\infty} \frac{i}{2}(\dot{s}_i+1)\textsc{h}
(m,\hat{m}, b\frac{\dot{b}_i}{\grave{b}_i})p_b=
\sum_c\sum_{i=1}^{\infty} \frac{i}{2}\textsc{h} (m,\hat{m},
a)\grave{p}_i\frac{\partial p_c}{\partial\dot{p}_i}$, where the sum
is given by all monomial $c$.

Make the changing $c=b\frac{\bar{b}_{i+j}}{\bar{b}_i \grave{b}_j}$
in the second summand. Then $p_c=p_b\frac{\bar{p}_{i+j}}{\bar{p}_i
\grave{p}_j}$ and $p_b=p_c\frac{\bar{p}_i
\grave{p}_j}{\bar{p}_{i+j}}= \bar{p}_i\grave{p}_j\frac{\partial
p_c}{\partial\bar{p}_{i+j}} (\bar{s}_{i+j}+1)^{-1}$. Thus
$\sum_b\sum_{i=1}^{\infty}\sum_{j=1}^{\infty}(\bar{s}_{i+j}+1)
\textsc{h}(m,\hat{m}, b\frac{\bar{b}_{i+j}}{\bar{b}_i
\grave{b}_j})p_b=\sum_c\sum_{i=1}^{\infty}\sum_{j=1}^{\infty}
\textsc{h}(m,\hat{m},c)\bar{p}_i\grave{p}_j\frac{\partial
p_c}{\partial\bar{p}_{i+j}}$.

Analogically,
$\sum_b\sum_{i=1}^{\infty}\sum_{j=1}^{\infty}(\grave{s}_{i+j}+1)
\textsc{h}(m,\hat{m}, b\frac{\grave{b}_{i+j}}{\grave{b}_i
\grave{b}_j})p_b=$ $\sum_c\sum_{i=1}^{\infty}\sum_{j=1}^{\infty}
\textsc{h}(m,\hat{m},c)\grave{p}_i\grave{p}_j\frac{\partial
p_c}{\partial\grave{p}_{i+j}}$ and
$\sum_b\sum_{i=1}^{\infty}\sum_{j=1}^{\infty}(\acute{s}_{i+j-1}+1)
\textsc{h}(m,\hat{m}, b\frac{\acute{b}_{i+j-1}}{\bar{b}_i
\bar{b}_j})p_b=$ $\sum_c\sum_{i=1}^{\infty}\sum_{j=1}^{\infty}
\textsc{h}(m,\hat{m},c)\bar{p}_i\bar{p}_j\frac{\partial
p_c}{\partial\acute{p}_{i+j-1}}$. Thus  $\frac{\partial
\acute{H}} {\partial\beta}= L_{\beta}H$.

$\Box$

Analogically, using lemma \ref{l3.2} instead lemma \ref{l3.1} we prove

\begin{theorem}\label{t4.2} $$\frac{\partial
\grave{H}}{\partial\gamma}=L_{\gamma}H,$$ где
$$L_{\gamma} =
\sum_{i=1}^{\infty}\dot{p}_i\frac{\partial}{\partial\grave{p}_i} + \sum_{ij}(
2\bar{p}_{i+j}\frac{\partial^2}{\partial\bar{p}_i\partial\grave{p}_j}+
2\grave{p}_{i+j}\frac{\partial^2}
{\partial\grave{p}_i\partial\grave{p}_j}+
\frac{1}{2}\acute{p}_{i+j-1}\frac{\partial^2}
{\partial\bar{p}_i\partial\bar{p}_j}).$$
\end{theorem}
The last theorem of the paper is

\begin{theorem}\label{t4.2} $$\frac{\partial H}{\partial\alpha}=
L_{\alpha}H,$$ where

$$L_{\alpha}= \sum_{ij}(\frac{i+j}{2}
\dot{p}_i\dot{p}_j\frac{\partial}{\partial\dot{p}_{i+j}}+
\frac{ij}{2}
\dot{p}_{i+j}\frac{\partial^2}{\partial\dot{p}_i\partial\dot{p}_j}+
|\bar{b}_j|
\dot{p}_i\bar{p}_j\frac{\partial}{\partial\bar{p}_{i+j}}+
i|\bar{b}_j|
\bar{p}_{i+j}\frac{\partial^2}{\partial\dot{p}_i\partial\bar{p}_j}+$$
$$2|\hat{b}_j|
\dot{p}_i\acute{p}_j\frac{\partial}{\partial\acute{p}_{i+j}}+
2i|\hat{b}_j|
\acute{p}_{i+j}\frac{\partial^2}{\partial\dot{p}_i\partial\acute{p}_j}+
2|\hat{b}_j|
\dot{p}_i\grave{p}_j\frac{\partial}{\partial\grave{p}_{i+j}}+
2i|\hat{b}_j|
\grave{p}_{i+j}\frac{\partial^2}{\partial\dot{p}_i\partial\grave{p}_j})+$$

$$\sum_{i+j=k+l}2(|\bar{b}_i\hat{b}_j\bar{b}_k\hat{b}_l|
\bar{p}_k\acute{p}_l\frac{\partial^2}{\partial\bar{p}_i\partial\acute{p}_j}+
|\bar{b}_i\hat{b}_j\bar{b}_k\hat{b}_l|
\bar{p}_k\grave{p}_l\frac{\partial^2}{\partial\bar{p}_i\partial\grave{p}_j})+$$

$$\sum_{i+j+1=k+l}
2(1+\delta_{kl})|\hat{b}_i\hat{b}_j\bar{b}_k\bar{b}_l|
\bar{p}_k\bar{p}_l\frac{\partial^2}{\partial\acute{p}_i\partial\grave{p}_j}+
\sum_{i+j=k+l+1}\frac{1}{2}(1+\delta_{ij})|\bar{b}_i\bar{b}_j\hat{b}_k\hat{b}_l|
\acute{p}_k\grave{p}_l\frac{\partial^2}{\partial\bar{p}_i\partial\bar{p}_j}+$$

$$\sum_{i+j=k+l}(1+\delta_{ij})(1+\delta_{kl})
(\frac{1}{4}|\bar{b}_i\bar{b}_j\bar{b}_k\bar{b}_l|
\bar{p}_k\bar{p}_l\frac{\partial^2}{\partial\bar{p}_i\partial\bar{p}_j}+
|\hat{b}_i\hat{b}_j\hat{b}_k\hat{b}_l|
\grave{p}_k\grave{p}_l\frac{\partial^2}{\partial\grave{p}_i\partial\grave{p}_j}+
|\hat{b}_i\hat{b}_j\hat{b}_k\hat{b}_l|
\acute{p}_k\acute{p}_l\frac{\partial^2}{\partial\acute{p}_i\partial\acute{p}_j})
$$

\end{theorem}

The proof use theorem \ref{t3.2} and is given by the same method
that the proof of theorem \ref{t4.1}.

\end{document}